\documentclass[12pt]{article}
\usepackage{amsmath,amsfonts, amsthm, mathrsfs}

\begin{document}

\addtolength{\textheight}{1pc}

\newcommand{\supp}{\operatorname{supp}}

\newtheorem{prop}{Proposition}[section]
\newtheorem{thm}{Theorem}[section]
\newtheorem{corr}{Corrolary}[section]
\newtheorem{lem}{Lemma}[section]
\newtheorem{defin}{Definition}[section]

% To get the enumeration characters right
\renewcommand{\theenumi}{(\roman{enumi})}
\renewcommand{\labelenumi}{\theenumi}

\makeatletter
\renewcommand{\theequation}{\arabic{section}.\arabic{equation}}
\@addtoreset{equation}{section}
\makeatother

% To get the theorem numbering right

%\bibliographystyle{plain}
%\nocite{*}

\begin{centering}
\null

\null

\LARGE On Hessian measures for non-commuting vector
fields.\normalsize\footnote{Research supported by Australian
Research Council grant}

\null

\large
Neil S Trudinger \\ Centre for Mathematics and its
Applications \\ Australian National University \\
\end{centering}
\normalsize

\begin{abstract}

Previous results on Hessian measures by Trudinger and Wang are
extended to the subelliptic case. Specifically we prove the weak
continuity of the 2-Hessian operator, with respect to local $L^1$
convergence, for a system of $m$ vector fields of step 2 and
derive gradient estimates for the corresponding $k$-convex
functions, $1 \le k \le m$.

\end{abstract}

\section{Introduction}

In the paper \cite{num14}, we introduced the notion of
$k$-convexity, $k = 1 \ldots n$, for functions $u$ defined on
domains $\Omega$ in Euclidean space, $\mathbb{R}^n$. Namely, for
$u \in C^2(\Omega)$, we call $u$ \emph{$k$-convex} in $\Omega$ if
\begin{equation}
F_j[u] := F_j(D^2u) := S_j(\lambda) \ge 0, \label{lab11}
\end{equation}
for $j = 1 \ldots k$, where $\lambda = \left( \lambda_1, \ldots,
\lambda_n \right)$ are the eigenvalues of the Hessian matrix
$D^2u$ of second derivatives of $u$ and $S_j$ denotes the $j$th
elementary symmetric function, that is
\begin{equation}
S_j(\lambda) = \sum_{i_1 < \ldots < i_j} \lambda_{i_1} \ldots
\lambda_{i_j}, \quad j = 1 \ldots n. \label{lab12}
\end{equation}
When there is no confusion we use the same notation $F_j$ for both
the operator and the function on $\mathbb{R}^n \times
\mathbb{R}^n$. Equivalently, $u$ is $k$-convex in $\Omega$ if $u$
is \emph{subharmonic} with respect to the operator $F_k$ and this
is the basis for our definition of $k$-convexity for non-smooth
functions in our sequel papers \cite{num15}, \cite{num16},
\cite{num17}. The core result in our paper \cite{num15}, is that
the mapping $u \mapsto F_k[u]$ is weakly continuous as a mapping
from $L_{loc}^1(\Omega)$ to $M_{loc}(\Omega)$, the space of
locally finite measures in $\Omega$, that is for any subdomain
$\Omega' \subset \subset \Omega$, $\eta \in C_0^0(\Omega')$ and
positive constant $\epsilon$, there exists a constant $\delta$
such that
\begin{equation}
\left| \int_{\Omega} \eta \left( F_k(u) - F_k(v) \right) \right| < \epsilon
\end{equation}
whenever
\begin{equation}
\int_{\Omega'} \left| u - v \right| < \delta, \quad \int_{\Omega'}
\left| u + v \right| < 1, \nonumber
\end{equation}
for arbitrary $k$-convex $u$ and $v$. This result enables us to
define for any locally integrable $k$-convex function $u$, the
\emph{Hessian measure}, $\mu_k[u]$, as an extension of $F_k[u]$.
In our first paper \cite{num14}, we only proved the continuity of
$\mu_k$ from $C^0\left( \Omega \right)$ to $M_{loc} \left( \Omega
\right)$ but this was enough for the cases $k > \frac{n}{2}$,
which included the Monge-Amp\'{e}re measure when $k = n$. The weak
continuity of the Monge-Amp\'{e}re measure is a fundamental result
of Aleksandrov (see eg. \cite{num12}).

In this paper we extend our results in \cite{num15} to the case of
non-commuting vector fields but only prove the corresponding weak
continuity for the case $k = 2$. Our approach follows
\cite{num15}, with some help from \cite{num16} and \cite{num17},
and is inspired by the recent paper \cite{num8} on the special
case of the Heisenberg group $\mathbb{H}^1$, by Guti\'{e}rrez and
Montanari, where the more restrictive approach in \cite{num14} was
adequate. To formulate the main theorem, we let $X = X_1 \ldots
X_m$ denote a system of vector fields in $\mathbb{R}^n$, that is
first order differential operators of the form
\begin{equation}
X_i = \sum_{j = 1}^n b^{ij} D_j
\end{equation}
with coefficients $b^{ij} \in C^{\infty} \left( \bar \Omega
\right)$ (although weaker regularity will suffice). Then, for $k =
1 \ldots m$, we call a function $u \in C^2 \left( \Omega \right)$,
$k$-convex, with respect to $X$ if
\begin{equation}
F_j[u] := F_j(X_s^2 u) := S_j(\lambda) \ge 0, \quad j = 1 \ldots k
\end{equation}
where now $\lambda = \left( \lambda_1 \ldots \lambda_m \right)$
denote the eigenvalues of the \emph{symmetric Hessian},
\begin{equation}
X_s^2 u = \left[ \frac{1}{2} (X_i X_j + X_j X_i ) u \right]_{i,j = 1 \ldots m}.
\end{equation}

Our hypotheses on the vector fields $X_1 \ldots X_m$ are that:
\begin{enumerate}
\item they are anti-self adjoint, namely \label{item1}
\begin{displaymath}
X_i^* = -X_i, \quad i = 1 \ldots m ;
\end{displaymath}
\item they satisfy the Hormander condition, namely the Lie algebra
   generated by them spans $\mathbb{R}^n$ and \label{item2}
\item the second commutators formed from any two vector fields
   vanish. \label{item3}
\end{enumerate}

These conditions will be automatically satisfied by the vector
fields generating an homogeneous group of Heisenberg type. We can
now state the main theorem.

\begin{thm} \label{thm11}
The mappings
\begin{equation}
u \mapsto F_2[u] + \alpha \sum_{i < j} [X_i, X_j]^2 u,
\end{equation}
for $u$ 2-convex in $\Omega$, are weakly continuous from
$L_{loc}^1 \left( \Omega \right)$ to $M_{loc} \left( \Omega
\right)$, for any constant $\alpha$.
\end{thm}

As mentioned above, the special case of the Heisenberg group
$\mathbb{H}^1$, given by
\begin{eqnarray}
X_1 u & = & D_1 u - \frac{1}{2} x_2 D_3 u, \nonumber \\
X_2 u & = & D_2 u + \frac{1}{2} x_1 D_3 u, \\
\left[ X_1 , X_2 \right] u & = & D_3 u \nonumber,
\end{eqnarray}
is proved in \cite{num8}. Here $L_{loc}^1 \left( \Omega \right)$
convergence is equivalent to local uniform convergence in $\Omega$
and the proof is much simpler.

Theorem \ref{thm11} enables us to assign a Borel measure
$\mu_2[u]$ to any $L^1_{loc}(\Omega)$ limit of smooth $2$-convex
functions, which extends $F_2[u]$ and is weakly continuous.
Letting $\Phi^2(\Omega)$ denote the space of such functions we
also see that the commutators $[ X_i, X_j ] u \in L_{loc}^2 (
\Omega )$ for $u \in \Phi^2 ( \Omega )$.

This paper is arranged as follows. In the next section we
generalize the basic divergence identity of Guti\'{e}rrez and
Montanari \cite{num8}, \cite{num9} on Heisenberg groups to vector
fields satisfying conditions \ref{item1} and \ref{item3}. A more
complete treatment in Carnot groups, with applications to
monotonicity, is given by Danielli, Garofalo, Nhieu and Tournier
in \cite{num5}; (see also \cite{num7}).  In Section 3, we employ
our approach in \cite{num15} to obtain integral estimates for the
subelliptic gradient $Xu$ under conditions \ref{item1} and
\ref{item2}, Theorem \ref{thm31}. As we cannot extend all of our
argument in \cite{num15}, we have to rely strongly on the
subelliptic potential estimates in \cite{num17}. In Section 4, we
carry out the arguments, again adapting \cite{num15} to the
non-commutative case, to conclude the local boundedness and weak
continuity of the functionals in Theorem~\ref{thm11}, thereby
completing the proof. Finally in Section 5, we extend our previous
results to the classes $\phi^k(\Omega)$ of $L_{loc}^1$ limits of
$k$-convex functions. Additional remarks at the ends of Sections 4
and 5 treat the removal of condition \ref{item3} and more general
definitions of $k$-convexity.

We are grateful for useful comments and discussions with N.
Chaudhuri, X-J Wang, A. Montanari and T. Nguyen.

\section{Divergence structure and monotonicity}

It is well known that if the vector fields $X_1 \ldots X_m$
commute, then
\begin{equation}
X_i F_k^{ij} (X^2 u) = 0, \quad j = 1 \ldots m, \label{lab21}
\end{equation}
where
\begin{equation}
X^2 u = X_s^2 u = [ X_i X_j u ]_{i,j = 1 \ldots m} \label{lab22}
\end{equation}
and
\begin{equation}
F_k^{ij}(r) := \frac{\partial}{\partial r_{ij}} F_k(r) \label{lab23}
\end{equation}
The identity (\ref{lab21}), which means that the columns of the
linearized coefficient matrix (\ref{lab23}) are divergence free,
was the basis for our approach in \cite{num14}, \cite{num15}. Now
suppose, more generally, that the second commutators formed from
any two vector fields vanish, that is for any $i, j = 1 \ldots m$,
\begin{eqnarray}
0 & = & [X_i, [X_i, X_j] ] \nonumber \\
  & = & X_i [X_i, X_j] - [X_i, X_j] X_i \nonumber \\
  & = & X_i X_i X_j - 2 X_i X_j X_i + X_j X_i X_i \nonumber \\
  & = & X_j(X_i X_i) + X_i ( X_i X_j - 2 X_j X_i ). \nonumber
\end{eqnarray}

Then, defining for any real matrix $r \in \mathbb{R}^n \times
\mathbb{R}^n$,
\begin{eqnarray}
\mathscr{F}_2(r) & := & \frac{1}{2} \left\{
    (r_{ii})^2 - r_{ij} r_{ji} + \frac{1}{2}
    (r_{ij} - r_{ji}) ^2 \right\} \label{lab24} \\
 & = & F_2[r] + \frac{3}{4} \sum_{i < j}
    (r_{ij} - r_{ji})^2, \nonumber
\end{eqnarray}
we have the identity
\begin{equation}
X_i \mathscr{F}_2^{ij}(X^2 u) = 0, \label{lab25}
\end{equation}
where
\begin{eqnarray}
\mathscr{F}_2^{ij}(r) & = & \frac{\partial \mathscr{F}_2}{\partial r_{ij}} (r) \label{lab26} \\
 & = & \text{(trace } r \text{)} \delta_{ij} + r_{ij} - 2 r_{ji}, \nonumber
\end{eqnarray}
which extends (\ref{lab21}) in the case $k = 2$. The identity
(\ref{lab25}) was discovered by Guti\'{e}rrez and Montanari
\cite{num8}, \cite{num9} for the Heisenberg groups $\mathbb{H}^n$,
(see also \cite{num5}, \cite{num7}). From (\ref{lab25}), we infer
the monotonicity formula for the operator $\mathscr{F}_2$ defined
by
\begin{equation}
\mathscr{F}_2[u] = \mathscr{F}_2(X^2 u), \label{lab27}
\end{equation}
extending Lemma 2.1 in \cite{num14} for $k = 2$.

\begin{lem} \label{lem21}
Let $u, v \in C^2(\Omega) \cap C^0(\bar \Omega)$ satisfy $u \le v$
in $\Omega$, $u = v$ on $\partial \Omega$ with the operator $F_2$
degenerate elliptic with respect to their sum $u + v$, that is
\begin{eqnarray}
\mathscr{F}_2^{ij}( X^2(u+v))\xi_i \xi_j & = & F_2^{ij}(X^2(u+v)) \xi_i \xi_j \label{lab28} \\
 & \ge & 0 \nonumber
\end{eqnarray}
for all $\xi \in \mathbb{R}^m$. Then, if the vector fields $X_1
\ldots X_m$ satisfy conditions \ref{item1} and \ref{item3}, we
have
\begin{equation}
\int_{\Omega} \mathscr{F}_2 [ v ] \le \int_{\Omega} \mathscr{F}_2 [ u ] \label{lab29}
\end{equation}
\end{lem}
\begin{proof}
By integration by parts and the identity (\ref{lab25}), we have,
for $u,v \in C^2( \bar \Omega)$,
\begin{eqnarray}
\int_{\Omega} \left( \mathscr{F}_2 [ u ] - \mathscr{F}_2 [v] \right)
 & = & \int_0^1 dt \int_{\Omega} \mathscr{F}_2^{ij}
    [ X^2 (tu + (1-t)v) ] X_iX_j(u-v) \nonumber \\
 & = & \int_0^1 dt \int_{\partial \Omega}
    \mathscr{F}_2^{ij} (X_i . \gamma ) X_j (u - v) \nonumber \\
 & = & \int_0^1 dt \int_{\partial \Omega}
    \mathscr{F}_2^{ij} (X_i . \gamma ) (X_j . \gamma)
    \left| D (u - v) \right| \nonumber \\
 & \ge & 0 \nonumber
\end{eqnarray}
Here $\gamma$ denotes the outer unit normal to $\partial \Omega$
and
\begin{displaymath}
X_i . \gamma = b^{ij} \gamma_j.
\end{displaymath}
The general case $u, v \in C^0 \left( \bar \Omega \right) \cap C^2
\left( \Omega \right)$ follows by approximation.
\end{proof}

More general version of Lemma~\ref{lem21} are presented in
\cite{num5}. For weak continuity with respect to $C^0 \left(
\Omega \right)$ and for groups of Heisenberg type we may proceed
exactly as in \cite{num14}. In the next section we present the
basic gradient estimates for $k$-convex functions needed to handle
the general case, for which we will not need Lemma~\ref{lem21}.

\section{Gradient Estimates}

In this section we provide the necessary gradient estimates for
our proof of weak continuity. For these we do not have to restrict
to the case $k = 2$ and moreover we only need to assume the vector
fields $X_1 \ldots X_m$ satisfy conditions \ref{item1} and
\ref{item2}. First we note that since $k$-convexity implies
1-convexity, $k$-convex functions $u$ are subharmonic with respect
to the sub-Laplacian associated with $X_1 \ldots X_n$, that is
\begin{equation}
\Delta_{X} u := X_i X_i u \ge 0 \label{lab31}
\end{equation}
in $\Omega$. From (\ref{lab31}) we infer immediately a bound from
above, namely, for any $\Omega' \subset \subset \Omega$,
\begin{equation}
\sup_{\Omega'} u \le C \int_{\Omega} \left| u \right| \label{lab32}
\end{equation}
where the constant $C$ depends on $X_1 \ldots X_m$ and
$\mathrm{dist}(\Omega', \partial \Omega)$. As we are only dealing
with local estimates in this paper, we will always assume, without
loss of generality, that $u \in L^1(\Omega)$. Following
\cite{num15}, our treatment of gradient estimates depends on the
relation between $k$-convexity and the subelliptic p-Laplacian
operators $\Delta_p$ defined by
\begin{equation}
\Delta_p u = X_i ( \left| X u \right| ^{p-2} X_i u ) \label{lab33}
\end{equation}
for $p > 1$.
\begin{lem} \label{lem31}
Let $u$ be $k$-convex in $\Omega$. Then $u$ is subharmonic with
respect to $\Delta_p$ for $p - 1 \le k \left(m-1\right) /
\left(m-k\right)$.
\end{lem}
\begin{proof}
Although this is just the special case $l=1$ in Lemma 4.2 of
\cite{num15}, we include it for completeness as it is simpler than
the cases $l>1$. We use the notation
\begin{equation}
S_{k,i} (\lambda) = S_k(\lambda) |_{\lambda_i = 0} \label{lab34}
\end{equation}
so that
\begin{equation}
S_{j,i} (\lambda) \ge 0 \label{lab35}
\end{equation}
for all $j \le k-1$, if $S_j(\lambda) \ge 0$, for all $j = 1,
\ldots k$, \cite{num15}. It follows then that
\begin{equation}
0 \le S_k(\lambda) = S_{k,i}(\lambda) +
  S_{k-1,i}(\lambda) \lambda_i,\quad i = 1, \ldots m, \label{lab36}
\end{equation}
whence
\begin{eqnarray}
-\lambda_i & \le & \frac{S_{k,i}}{S_{k-1,i}} (\lambda) \label{lab37} \\
 & \le & \frac{(m-k)}{k(m-1)} S_{1,i}(\lambda) \nonumber
\end{eqnarray}
by MacLaurins' inequality for ratios of elementary symmetric
functions. Consequently, if $u$ is $k$-convex,
\begin{eqnarray}
\Delta_p u & = & X_i ( \left| Xu \right|^{p-2} X_i u ) \label{lab38} \\
 & = & \left| Xu \right|^{p-2} \left\{ \Delta_Xu + (p-2)
    \frac{X_i u X_j u}{\left| Xu \right|^2} X_i X_j u \right\} \nonumber \\
 & \ge & \left| Xu\right|^{p-2}\left\{\Delta_Xu+(p-2)
    \,\lambda_{\min}\!\left(X_s^2u\right)\right\}\nonumber\\
 & \ge & 0, \nonumber
\end{eqnarray}
for $p-1 \le k \left(m-1\right) / \left(m-k \right)$, by taking
$\lambda_1 \ldots \lambda_m$ in (\ref{lab37}) to be the
eigenvalues of $X^2_s u$.
\end{proof}
Note that Lemma~\ref{lem31} also includes the case $k = m, p =
\infty$, when $\Delta_p$ is the subelliptic $\infty$-Laplacian
\begin{equation}
\Delta_{\infty} u = X_i u X_j u X_{ij} u \label{lab39}
\end{equation}
Our gradient estimates now follow immediately from \cite{num17}
but to express them we need the concept of homogeneous dimension.
For our purposes here, we define the C-C (Carnot-Caratheodory)
metric induced from the vector fields $X_1 \ldots X_m$ by
\begin{equation}
d(x,y) = \inf \{ T > 0 | \exists \text{ a sub-unitary }
  \gamma: [0,T] \rightarrow \mathbb{R}^n \nonumber
\end{equation}
\begin{equation}
\quad \quad \text{ with } \gamma(0) = x, \gamma(T) = y \}, \label{lab310}
\end{equation}
where a piecewise $C^1$ curve $\gamma : [0,T] \rightarrow
\mathbb{R}^n$ is said to be sub-unitary, with respect to $X_1
\ldots X_m$, if for every $\xi \in \mathbb{R}^n$ and $t \in (0,
T)$,
\begin{equation}
\left| \gamma ' (t) . \xi \right| ^2 \le \sum_{i = 1}^m
  (X_i (\gamma(t)). \xi)^2 \label{lab311}.
\end{equation}

Let $B_R(x)$ denote the C-C ball $\left\{y \in \mathbb{R}^n |
d(x,y) < R \right\}$, and let $\Omega$ be a bounded domain in
$\mathbb{R}^n$. The the fundamental result of Nagel, Stein and
Wainger \cite{num11} asserts that there exist positive constants
$C, R_0$ and positive integer $Q$, depending on $X$ and $\Omega$
such that
\begin{equation}
\left| B_{tR} (x) \right| \ge C t^Q \left| B_R(x) \right| \label{lab312}
\end{equation}
for any $x \in \Omega$, $t \in (0,1)$ and $R < R_0$, where $| |$
denotes the Lebesgue volume. The number $Q$ ($\ge n$), is chosen
as the least integer for which (\ref{lab312}) holds and is called
the \emph{homogeneous dimension} of $X$ in $\Omega$. For
(\ref{lab312}) we only need the Hormander condition \ref{item2}
and we could replace it more generally in this paper by simply the
validity of (\ref{lab312}).

\begin{thm} \label{thm31}
For any $k$-convex function $u$ in $\Omega$, and subdomain
$\Omega' \subset \subset \Omega$, we have the estimates
\begin{equation}
|| X u||_{L^q(\Omega')} \le C \left( \int_{\Omega}
   \left| u \right| \right) \label{lab313}
\end{equation}
\begin{equation}
\int_{\Omega'} \left| Xu \right|^r \Delta_X u \le C \left(
\int_{\Omega} \left| u \right| \right)^{1+r} \label{lab314}
\end{equation}
for $1 \le q < Q k (m-1)/(Q-1)(m-k)$, $0 \le r < m (k-1)/(m-k)$
where $C$ depends on $\Omega$, $\Omega'$, $X_1 \ldots X_m$ and $q$
or $r$ as appropriate.
\end{thm}
\begin{proof}
The estimate (\ref{lab313}) follows from Lemmas~\ref{lem31} and
\cite{num17}, Lemma~3.9. For (\ref{lab314}), we have from
Lemma~\ref{lem31},
\begin{equation}
\left| Xu \right|^r \Delta_X u \le
 \frac{m(k-1)}{m(k-1) - r(m-k)} \Delta_p u \label{lab315}
\end{equation}
and since
\begin{equation}
\int_{\Omega} \eta \Delta_p u \le \int_{\Omega}
 \left| X \eta \right| \left| X u \right| ^{p-1}, \label{lab316}
\end{equation}
for any $\eta \ge 0, \in C_0^1(\Omega)$, we infer (\ref{lab314})
from (\ref{lab313}).
\end{proof}
By using the subelliptic Sobolev inequality \cite{num2},
\cite{new11}, we obtain corresponding $L^p$ estimates, namely
\begin{equation}
|| u ||_{L^p(\Omega')} \le C || u ||_{L^1(\Omega)} \label{lab317}
\end{equation}
where
\begin{displaymath}
1 \le p < \frac{Q k (m-1)}{(Q - 1)m - (Q + m - 2)k}
\end{displaymath}
for $(Q - 1)m \ge (Q + m - 2) k$ and $p = \infty$ if $(Q - 1)m <
(Q + m - 2)k$. For this last case we have a H\"{o}lder estimate
\cite{new12}, \cite{num17},
\begin{equation}
\sup_{\Omega'} \frac{\left|u(x) - u(y)\right|}{\left| d(x,y)
\right|^{\alpha}} \le C ||u||_{L^1(\Omega)} \label{lab318}
\end{equation}
where
\begin{displaymath}
\alpha = \frac{k(Q+m-2) - m(Q-1)}{k(m-1)},
\end{displaymath}
if $k < m$. When $k=m$, we may take $q=\infty$ in (\ref{lab318}),
$\alpha=1$ in (\ref{lab313}) as the estimates (\ref{lab318}) are
uniform in $q<\infty$.

\section{Weak Continuity}

In this section, we complete the proof of Theorem~\ref{thm11}.
First we prove a local bound for $\mathscr{F}_2$. For convenience
we use the notation
\begin{equation}
\mathscr{E}_2 [ u ] = \sum_{i < j} \left( [ X_i, X_j] u \right)^2, \label{ref41}
\end{equation}
so that
\begin{equation}
\mathscr{F}_2 [ u ] = F_2 [ u ] + \frac{3}{4} \mathscr{E}_2 [ u ]
\end{equation}
\begin{lem} \label{lem41}
Let $u \in C^2 ( \Omega )$ be 2-convex in $\Omega$ with respect to
$X_1, \ldots X_m$, satisfying hypotheses \ref{item1} to
\ref{item3}. Then, for any subdomain $\Omega' \subset \subset
\Omega$, we have
\begin{equation}
\int_{\Omega'} \mathscr{F}_2 [ u ] \le C
 \left( \int_{\Omega} \left| u \right| \right) ^2, \label{lab43}
\end{equation}
where $C$ depends on $dist(\Omega', \partial \Omega)$, $X_1,
\ldots X_m$.
\end{lem}

In particular, Lemma~\ref{lem41} provides a local $L^2$ estimate
for the commutators $[X_i, X_j]u$.

\begin{proof}
Letting $\eta > 0 \in C_0^1 ( \Omega )$ be a cut-off function, we
have
\begin{eqnarray}
\int_{\Omega} \eta^2 \mathscr{F}_2 [ u ]
 & = & \frac{1}{2} \int_{\Omega} \eta^2
   \mathscr{F}_2^{ij} X_i X_j u \label{lab44} \\
 & = & - \frac{1}{2} \int_{\Omega}
   \mathscr{F}_2^{ij} X_i \eta^2 X_j u \nonumber \\
 & = & - \frac{1}{2} \int_{\Omega} \left\{ F_2^{ij} X_i \eta^2 X_j u
    + \frac{3}{2} [X_i, X_j] u X_i \eta^2 X_j u \right\} \nonumber \\
 & \le & C \int_{\Omega} \left\{ \left| X \eta^2 \right|
   \left| X u \right| \Delta_X u + \eta \mathscr{E}_2^{\frac{1}{2}}
    [u] \left| X \eta \right| \left| X u \right| \right\}. \nonumber
\end{eqnarray}
Consequently,
\begin{equation}
\int_{\Omega} \eta^2 \mathscr{F}_2 [u] \le C \int_{\Omega} \left\{
\left| X \eta^2 \right| \left| Xu \right| \Delta_X u + \left| X
\eta \right|^2 \left| X u \right|^2 \right\} \label{lab45}
\end{equation}
and (\ref{lab43}) follows from Theorem~\ref{thm31}.
\end{proof}

\begin{proof}[Proof of Theorem~\ref{thm11}]

Letting $u$ and $v$ be 2-convex in $\Omega$, $u_t = tu + (1 - t)v,
0 \le t \le 1$ and $\eta \ge 0, \in C_0^2(\Omega)$, we now have
\begin{eqnarray}
 & & \int_{\Omega} \eta \left( \mathscr{F}_2 [u] - \mathscr{F}_2[v] \right) \label{lab46} \\
 & = & \int_0^1 dt \int_{\Omega} \eta \mathscr{F}_2^{ij} [u_t] X_i X_j (u - v) \nonumber \\
 & = & \frac{1}{2} \int_{\Omega} \eta \mathscr{F}_2^{ij} [u + v] X_i X_j (u - v) \nonumber \\
 & = & - \frac{1}{2} \int_{\Omega} \mathscr{F}_2^{ij} X_i \eta X_j (u - v) \nonumber \\
 & = & - \frac{1}{2} \int_{\Omega} \left\{ \mathscr{F}_2^{ji} + 3[X_i, X_j](u+v) \right\} X_i \eta X_j (u - v) \nonumber \\
 & = & \frac{1}{2} \int_{\Omega} \left\{ \mathscr{F}_2^{ij} X_i X_j \eta (u - v) - 3[X_i, X_j] (u + v) X_i \eta X_j (u - v) \right\} \nonumber \\
 & = & \frac{1}{2} \int_{\Omega} F_2^{ij} [u + v] (X_i X_j \eta) (u - v) \nonumber \\
 &   & \quad + \frac{3}{4} \int_{\Omega} [X_i, X_j] (u + v) \left\{(X_i X_j \eta)(u - v) - 2 X_i \eta X_j (u - v) \right\} \nonumber \\
 & := & \frac{1}{2} I_1 + \frac{3}{4} I_2 \nonumber
\end{eqnarray}
The estimation of the integral $I_1$ is similar to the
corresponding term in the commuting case \cite{num15}. Namely
\begin{eqnarray}
\left| I_1 \right| & = & \left| \int_{\Omega} F_2^{ij} [u + v] (X_i X_j \eta) (u - v) \label{lab47} \right| \\
 & \le & \int_{\Omega} \Delta_X (u + v) \left| X^2 \eta \right| \left| u - v \right|. \nonumber
\end{eqnarray}
Denoting
\begin{equation}
\delta = \int_{\Omega} \left| u - v \right|, \quad K = \int_{\Omega} \left| u + v \right|, \label{lab48}
\end{equation}
we have, for any $\epsilon > 0$,
\begin{displaymath}
\left| u - v \right| < \epsilon \nonumber
\end{displaymath}
except on a set $A_{\epsilon}$ of measure $\left| A_{\epsilon}
\right| \le \delta / \epsilon$. We then estimate for a further
cut-off function $\tilde{\eta} \in C_0^1(\Omega)$, $0 \le
\tilde{\eta} \le 1$,
\begin{eqnarray}
 & & \int_{\Omega} \tilde{\eta} \Delta_X (u + v) (u - v)^+ \label{lab49} \\
 & \le & \epsilon \int_{\Omega} \tilde{\eta} \Delta_X (u + v) + \int_{\Omega} \tilde{\eta} (u - v - \epsilon)^+ \Delta_X (u + v) \nonumber \\
 & = & \epsilon \int_{\Omega} \tilde{\eta} \Delta_X (u + v) - \int_{A_{\epsilon}} X_i (u + v) X_i \left\{ \tilde{\eta} (u - v - \epsilon)^+ \right\} \nonumber \\
 & \le & C \left( \epsilon + \left| A_{\epsilon} \right| ^{1 - \frac{2}{q}} \right) \nonumber \\
 & \le & C \left\{ \epsilon + \left( \frac{\delta}{\epsilon} \right) ^{1 - \frac{2}{q}} \right\} \nonumber
\end{eqnarray}
by Theorem~\ref{thm31}, where $q$ is chosen so that
\begin{displaymath}
2 < q < \frac{2 Q (m - 1)}{(Q - 1)(m - 2)}
\end{displaymath}
and $C$ depends on $\tilde{\eta}$, $K$, $\delta$, $X_1, \ldots
X_m$. It then follows that
\begin{equation}
\left| I_1 \right| \le C \left\{ \epsilon + \left( \frac{\delta}{\epsilon} \right)^{1 - \frac{2}{q}} \right\} \label{lab410}
\end{equation}
where $C$ depends on $\eta$, $K$, $\delta$, $X_1, \ldots X_m$. To
estimate $I_2$, we first use the bound for $\mathscr{E}_2$ in
Lemma~\ref{lem41}, to obtain
\begin{equation}
I_2 \le C \left\{ \int_{\supp \eta} \left| X(u - v) \right|^2 +
\left| u - v \right|^2 \right\}^\frac{1}{2}, \label{lab411}
\end{equation}
where $C$ depends on $\eta$, $K$, $\delta$, $X_1, \ldots X_m$. The
first part of the above integral may be estimated similarly to
$I_1$, since
\begin{eqnarray}
\int_{\Omega} \tilde{\eta} \left| X (u-v) \right|^2 & = & - \int_{\Omega} X_i (\tilde{\eta} X_i (u - v) ) (u - v) \label{lab412} \\
 & \le & \int_{\Omega} \left\{ \left|X \tilde{\eta} \right| \left| X (u-v) \right| + \tilde{\eta} \Delta_X (u-v) \right\} \left| u - v \right| \nonumber \\
 & \le & \int_{\Omega} \tilde{\eta} \Delta_X (u - v) \left| u - v \right| + C \left\{ \int_{\supp \tilde{\eta}} \left| u - v \right|^2 \right\}^\frac{1}{2}, \nonumber
\end{eqnarray}
while the second parts of (\ref{lab411}) and (\ref{lab412}) are
handled readily by the estimate (\ref{lab317}). It follows that
\begin{equation}
I_2 \le C \left\{ \epsilon + \left( \frac{\delta}{\epsilon}
\right) ^{1 - \frac{2}{q}} \right\}^{\frac{1}{2}} \label{lab413}.
\end{equation}
With appropriate choice of $\delta$, we then conclude
Theorem~\ref{thm11} from the estimates (\ref{lab410}) and
(\ref{lab413}).
\end{proof}
\textbf{Remark.} By inspection of the above proofs, we see that
condition \ref{item3} may be weakened to only requiring that the
vector fields
$$
Y_j := \sum_{i = 1}^m [X_i, [X_i, X_j] ], \quad j = 1 \ldots m,
$$
lie in the span of $X_i, [X_i, X_j]$, $i,j = 1 \ldots m.$ Without
condition \ref{item3} additional terms
\begin{equation}
- \frac{1}{2} \int_{\Omega} \eta^2 Xu.Yu \label{lab414}
\end{equation}
and
\begin{equation}
- \frac{1}{2} \int_{\Omega} \eta X(u-v).Y(u+v) + \frac{1}{2}
\int_\Omega (u-v)X \eta . Y(u+v) \label{lab415}
\end{equation}
will arise in the right hand sides of (\ref{lab44}) and
(\ref{lab46}) respectively and these are then automatically
controlled. More general hypotheses are clearly possible. Note
also that if $Y$ commutes with $X$, then the quantity
\begin{equation}
\mathscr{F}_{2}^{*} [u] := \mathscr{F}_2[u] + \frac{1}{2} Xu.Yu \label{lab416}
\end{equation}
also satisfies the monotonicity property (\ref{lab29}). For the
Engel group, this was first observed in \cite{num7}. More
generally if we integrate by parts in (\ref{lab415}), we find
\begin{eqnarray}
 &   & \int_\Omega \eta X(u-v).Y(u+v) \label{lab417} \\
 & = & - \int_\Omega \eta(u+v) Y.X(u-v) + (u+v) Y \eta .X(u-v) \nonumber \\
 & = & \int_\Omega \eta \{X(u+v).Y(u-v) + (u+v) Z(u-v) \} \nonumber \\
 &   & \quad - \int_\Omega (u+v) \{ Y \eta . X (u-v) - X \eta.Y(u-v) \} \nonumber
\end{eqnarray}
where $Z$ is the vector field given by
\begin{equation}
Z = X.Y - Y.X=0. \label{lab418}
\end{equation}
Inserting (\ref{lab417}) into (\ref{lab46}) and (\ref{lab415}) we
may infer weak continuity results, in the absence of condition
\ref{item3}, with respect to stronger topologies. For example
$F_2$ will be weakly continuous with respect to
$L^{1}_{loc}(\Omega)$ on 2-convex functions, with uniformly
bounded $X$- and $Y$-gradients in $L^\infty_{loc}(\Omega)$ and
$L^1_{loc}(\Omega)$ respectively.

Finally without hypothesis \ref{item3} in Lemma~2.1, we obtain,
from the proof,
\begin{eqnarray}
\int_\Omega (\mathscr{F}_2[u] - \mathscr{F}_2[v] ) & \ge & - \frac{1}{2} \int_\Omega Y(u+v).X(u-v) \nonumber \\
 & = & - \frac{1}{4} \int_\Omega \{Y(u+v).X(u-v) + Y(u-v).X(u+v)\} \nonumber
 \end{eqnarray}
so that in general
\begin{equation}
\int_\Omega ( \mathscr{F}_2^*[u] - \mathscr{F}_2^*[v] ) \ge 0.
\end{equation}
See \cite{num5} for a more thorough analysis of monotonicity. In a
similar fashon, we may combine (\ref{lab415}) and (\ref{lab417})
in (\ref{lab46}) to conclude,
\begin{eqnarray}
\int_\Omega \eta ( \mathscr{F}_2^*[u] - \mathscr{F}_2^*[v] )
 & = & \frac{1}{2} I_1 + \frac{3}{4} I_2 + \frac{3}{4} \int_\Omega (u-v)X\eta.Y(u+v)\qquad \\
 & &{}-  \frac{1}{4} \int_\Omega (u-v) X(u+v).Y \eta\,.\nonumber
\end{eqnarray}

\section{General $k$-convex functions}

For our purposes here, we define a function $u \in
L_{loc}^1(\Omega)$ to be $k$-convex in $\Omega$, with respect to
the system of vector fields $X = \left( X_1 \ldots X_m \right)$
if, for any $\Omega' \subset \subset \Omega$, there exists a
sequence of $k$-convex functions in $C^2 \left( \Omega' \right)$
converging to $u$ in $L^1_{loc} \left( \Omega' \right)$. We
designate the general class of $k$-convex functions in $\Omega$ by
$\phi^k \left( \Omega \right)$. The estimates of the preceeding
sections then extend as regularity properties and imbeddings of
$\phi^k \left( \Omega \right)$. In particular if the system $X$
satisfies conditions \ref{item1} and \ref{item2}, we have from
Theorem~\ref{thm31} that the distributional derivatives $Xu \in
L_{loc}^q \left( \Omega \right)$ for
\begin{equation}
1 \le q < \frac{Q k (m - 1)}{(Q - 1) (m - k)} \label{lab51}
\end{equation}
and that $\phi^k ( \Omega )$ imbeds continuously in the Sobolev space,
\begin{equation}
S_{loc}^{1,q} ( \Omega ) = \left\{ u \in L_{loc}^1 \left( \Omega
\right) | Xu \in L^q_{loc} ( \Omega ) \right\}, \label{lab52}
\end{equation}
as defined, for example, in \cite{num2}, \cite{num17}. Moreover for
\begin{equation}
k > \frac{(Q - 1) m}{(Q + m - 2)}, \label{lab53}
\end{equation}
we see from (\ref{lab318}) that $\phi^k (\Omega)$ imbeds
continuously in the H\"{o}lder space $C^{0,\alpha}(\Omega)$ where
\begin{equation}
\alpha = \frac{k(Q+m-2) - m(Q-1)}{k(m-1)}, \label{lab54}
\end{equation}
if $k < m$ ($\alpha < 1$ if $k = m$). It follows from \cite{num10}
that we can also take $\alpha = 1$ when $k = m$ and $X$ generates
the Lie algebra of a Carnot group.

For $k \ge 2$, the symmetric Hessian $X_s^2 u$ consists of signed
Radon measures. This follows exactly as in the Euclidean case
\cite{num4}, \cite{num15} but it may also be observed directly
from the degenerate ellipticity of $F_2$,
\begin{equation}
\Delta_X u I - X_s^2 u \ge 0. \label{lab55}
\end{equation}
By taking limits, we see that the above quantities are measures if
$u \in \phi^k (\Omega)$. If we also assume $X$ satisfies
\ref{item3}, we obtain from Lemma~\ref{lem41}, that the
commutators $[X_i, X_j] u \in L^2_{loc} (\Omega)$ for $u \in
\phi^k(\Omega)$, $k \ge 2$. Hence the full Hessian $X^2 u$
comprises Radon measures. From Theorem~\ref{thm31}, we also infer
that $\phi^k(\Omega)$, $k \ge 2$, imbeds continuously into the
Euclidean Sobolev space $W_{loc}^{1,2} (\Omega)$ if condition
\ref{item3} is strengthened to all second order commutators
vanishing. Also if $X$ generates the Lie algebra of a Carnot group
as in \cite{num1}, then from the H\"{o}lder estimate
(\ref{lab318}) and the weak differentiability result in
\cite{num1} we conclude that functions in $\phi^k(\Omega)$ will be
twice differentiable almost everywhere in $\Omega$, with respect
to the system $X$, if $k$ satisfies (\ref{lab53}). This extends
the corresponding result for the convex case $k = m$ in
\cite{num5}, \cite{num9}, \cite{num10} and the Euclidean case $k >
n/2$ in \cite{num4}.

Finally we may define the Hessian measure $\mu_2[u]$, with respect
to $X$, for any function $u \in \phi^2(\Omega)$ by
\begin{equation}
\int_{\Omega} \eta d \mu_2[u] = \lim_{m \rightarrow \infty} \int_{\Omega} \eta d \mu_2 [u_m]
\end{equation}
where $\eta \in C_0^0(\Omega)$, $\supp  \eta \subset \Omega'
\subset \subset \Omega$ and $\{u_m\} \subset C^2(\Omega')$ is a
sequence of $k$-convex functions converging to $u$ in
$L^1_{loc}(\Omega')$. By Theorem~\ref{thm11} $\mu_2$ is well
defined and weakly continuous with respect to convergence in
$L_{loc}^1(\Omega)$, that is if $\{u_m\} \subset \phi^2(\Omega)$,
converges to $u \in \phi^2(\Omega)$ in $L_{loc}^1(\Omega)$, the
corresponding sequence of measures $\{ \mu_2[u_m]\}$ converges
weakly to $\mu_2[u]$. Note that the case $\alpha = \frac{1}{2}$ in
Theorem~\ref{thm11} also shows that the sum of the principal $2
\times 2$ minors of the full Hessian $X^2 u$ also extends as a
weakly continuous measure on $\phi^2(\Omega)$. From
Lemma~\ref{lem21}, we also conclude a more general monotonicity
property, namely that if $u,v \in \phi^2(\Omega)$ satisfy $u \le
v$ in $\Omega$, $u = v$ continuously on $\partial \Omega$, then
\begin{equation}
\mu_2[v](\Omega) + \frac{3}{4} \mathscr{E}_2[v] \le
\mu_2[u](\Omega) + \frac{3}{4} \mathscr{E}_2[u].
\end{equation}

\subsection*{General subharmonic functions}

More generally we may define subharmonic functions along the lines
of \cite{num15}, \cite{num16}, \cite{num17}. In particular we
define an upper-semicontinuous function $u: \Omega \mapsto [
-\infty, \infty )$ to be subharmonic with respect to the operator
$F_k$ if $u$ satisfies $F_k[u] \ge 0$ in the viscosity sense, that
is for any quadratic polynomial $q$ for which the difference $u -
q$ has a finite local maximum at a point $y \in \Omega$, we have
$F_k[q](y) \ge 0$. For smooth vector fields and $k = 1$, this is
also equivalent to our definition in \cite{num17}, which
corresponds to the traditional definition of subharmonicity. A
$k$-convex function, as defined above by approximation, will be
equivalent to a subharmonic function and moreover the estimates of
Section~3 extend to the class of proper subharmonic functions. For
Carnot groups of step 2, it follows from \cite{num10} that proper
subharmonic functions will also be $k$-convex but we would expect
this characterization to hold more generally. The equivalence of
various definitions in the convex case, $k=m$, for Carnot groups
is treated in the papers \cite{new4}, \cite{new10}, \cite{new13},
\cite{num10}.

%\bibliography{OnHessianMeasures}

\begin{thebibliography}{10}

\bibitem{num1}
L.~Ambrosio and V.~Magnani. Weak differentiability of {BV}
functions on stratified groups.
 {\em Math. Zeit.}, (245):123 -- 153, 2003.

\bibitem{num2}
L.~Capogna, D.~Danielli, and N.~Garofalo.
 An embedding theorem and the {H}arnack inequality for nonlinear
  subelliptic equations.
 {\em Comm. Part. Diff. Eqns}, (18):1765 -- 1794, 1993.

\bibitem{num4}
N.~Chaudhuri and N.S. Trudinger.
 An {A}lexsandrov type theorem for $k$-convex functions.
 {\em Bull. Aust. Math. Soc.}, 71:305--314, 2005.

\bibitem{new4}
D.~Danielli, N.~Garofalo, and D.M.~Nhieu.
 Notions of convexity in Carnot groups.
 {\em Comm. Anal. Geom.}, 11(2):263--341, 2003.

\bibitem{num5}
D.~Danielli, N.~Garofalo, D.M. Nhieu, and F.~Tournier.
 The theorem of {B}usemann-{F}eller-{A}lexandrov in {C}arnot groups.
 {\em Comm. Anal. Geom}, 12(4):853--886, 2004.

\bibitem{num6}
G.B. Folland and E.M. Stein.
 Hardy spaces on homogeneous groups.
 {\em Princeton Univ. Press}, 1982.

\bibitem{num7}
N.~Garofalo and F.~Tournier.
 New properties of convex functions in the {H}eisenberg group.
 {\em Trans. Amer. Math. Soc.}, to appear.

\bibitem{num8}
C.E. Guti\'{e}rrez and A.~Montanari.
 Maximum and comparison principles for convex functions on the
  {H}eisenberg group.
 {\em Comm. Part. Diff. Eqns.}, 29:1305--1334, 2004.

\bibitem{num9}
C.E. Guti\'{e}rrez and A.~Montanari.
 On the second order derivatives of convex functions on the
  {H}eisenberg group.
 {\em Ann. Scuola Norm. Sup. Pisa Cl. Sci (5)}, 3(2):349--366, 2004.

\bibitem{new10}
P.~Juutinen, G.~Lu, J.J. Manfredi, and B.~Stroffolini.
 Convex functions on Carnot groups, (preprint).

\bibitem{new11}
G.~Lu.
 Weighted Poncar\'e and Sobolev inequalities for vector fields
 statisfying Hormander's condition and applications.
 {\em Rev. Mat. Iberoamericana} 8:367--439, 1992.

\bibitem{new12}
G.~Lu.
 Embedding theorems into Lipschitz and BMO spaces and applications
 to quasilinear subelliptic differential equations.
 {\em Publ. Mat.} 40:301--329, 1996.

\bibitem{new13}
G.~Lu, J.J. Manfredi and B.~Stroffolini.
 Convex functions on the Heisenberg Group.
 {\em Calc. Var. Partial Differential Equations} 19:1--22, 2004.

\bibitem{num10}
V.~Magnani.
 Lipschitz continuity, {A}leksandrov theorem and characterizations for
  {H}-convex functions.
 {\em Math. Annalen} (to appear).

\bibitem{num11}
A.~Nagel, E.M. Stein, and S.~Wainger.
 Balls and metrics defined by vector fields.
 {\em Acta. Math.}, 155:103--147, 1985.

\bibitem{num12}
A.V. Pogorelov.
 Monge-{A}mp\'{e}re equations of elliptic type.
 {\em Noordhoff, Groningen}, 1964.

\bibitem{num13}
E.M. Stein.
 {\em Harmonic {A}nalysis}.
 Princeton Univ Press, 1993.

\bibitem{num14}
N.S. Trudinger and X.J. Wang.
 Hessian {M}easures {I}.
 {\em Topol. Methods Nonlinear Anal.}, 10:225--239, 1997.

\bibitem{num15}
N.S. Trudinger and X.J. Wang.
 Hessian {M}easures {II}.
 {\em Ann. Math}, 150:579--604, 1999.

\bibitem{num16}
N.S. Trudinger and X.J. Wang.
 Hessian {M}easures {III}.
 {\em J. Funct. Anal.}, 193:1--23, 2002.

\bibitem{num17}
N.S. Trudinger and X.J. Wang.
 On the weak continuity of elliptic operators and applications to
  potential theory.
 {\em Amer. J. Math}, 124:369--410, 2002.

\end{thebibliography}

Centre for Mathematics and its Applications, Australian National University, Canberra, ACT 0200, Australia. Email: neil.trudinger@anu.edu.au
\end{document}